\documentclass[11pt]{article}
\usepackage{german}
\usepackage[T1]{fontenc}
\usepackage[ansinew]{inputenc}
\usepackage{latexsym}
\usepackage{textcomp}
\usepackage{color}
\usepackage{amscd}
\usepackage{amsmath}
\input{amssym.def}
\input{amssym}

     \newcommand{\fa}{\goth{a}}
     \newcommand{\fb}{\goth{b}}

     \newcommand{\ff}{\goth{f}}

     \newcommand{\fo}{\goth{o}}
     \newcommand{\fp}{\goth{p}}

     \newcommand{\fx}{\goth{x}}

     \newcommand{\fB}{\goth{B}}
     
     \newcommand{\fD}{\goth{D}}

     \newcommand{\fG}{\goth{G}}
     \newcommand{\fH}{\goth{H}}

     \newcommand{\C}{\Bbb{C}}

     \newcommand{\Q}{\Bbb{Q}}
     \newcommand{\R}{\Bbb{R}}

     \newcommand{\Z}{\Bbb{Z}}

    \newcommand{\ol}[1]{\overline{#1}}

    \newcommand{\ti}[1]{\tilde{#1}}
    \newcommand{\group}[1]{\langle{#1}\rangle}

    \newcommand{\ot}{\otimes}
    \newcommand{\me}{^{-1}}
    \newcommand{\mal}{^{\times}}
    
    \newcommand{\df}{\stackrel{\mathrm{def}}{=}}
    \newcommand{\mr}{\mathrm}
    
    \newcommand{\ilim}[1]{\raisebox{-1mm}{$\lim\atop{\leftarrow\atop\raisebox{0.5mm}{{\tiny $#1$}}}$}}

    \newcommand{\zp}{{\Bbb{Z}_p}}
    \newcommand{\zl}{{\Bbb{Z}_l}}
    \newcommand{\qp}{{\Bbb{Q}_p}}

    \newcommand{\onto}{\twoheadrightarrow}
    \newcommand{\lto}{\longrightarrow}
    \newcommand{\da}{\downarrow}

    \newcommand{\pl}{\parallel}
    \newcommand{\pht}{\phantom}

    \def\daz#1{#1\da\pht{#1}}

    \newcommand{\sda}{\da_{\raisebox{-1.4mm}{$\hsp{-1}\check{}$}}}

    \newcommand{\mapsdto}{\raisebox{1.75mm}{-}\hsp{-2.55}\da}

    \newcommand{\ga}{\gamma}
    \newcommand{\Ga}{\Gamma}
    \newcommand{\si}{\sigma}
    \newcommand{\Si}{\Sigma}
    \newcommand{\de}{\delta}
    \newcommand{\De}{\Delta}
    \newcommand{\ze}{\zeta}
    \newcommand{\la}{\lambda}
    
    \newcommand{\be}{\beta}
    \newcommand{\Om}{\Omega}
    
    \newcommand{\al}{\alpha}
    \newcommand{\ve}{\varepsilon}

    \newcommand{\noi}{\par\noindent}
    \newcommand{\sn}{\par\smallskip\noindent}
    \newcommand{\mn}{\par\medskip\noindent}
    \newcommand{\bn}{\par\bigskip\noindent}
    \newcommand{\bbn}{\par\bigskip\bigskip\noindent}
    
    \newcommand{\Section}[2]{\bbn {\large #1\,. \ {\sc #2}}
                             \nopagebreak
                             \nz}

    \newcommand{\nf}[2]{\\[1.5ex]
                        \bmp{1cm}
                         (#1)
                        \emp 
                        \bmp{13.5cm}
                         \bct
                          $#2$
                         \ect
                        \emp\\[1.5ex]
         }

    \newcommand{\sss}{\scriptstyle}
    \newcommand{\nz}{\\[1ex]}

    \newcommand{\hsp}[1]{\hspace*{#1mm}}

    \newcommand{\mmargin}{
     \textheight 230truemm
     \textwidth 155truemm
     \topmargin -10truemm
     \oddsidemargin 5truemm
     \evensidemargin 5truemm
     }

    \newcommand{\bmp}{\begin{minipage}}
    \newcommand{\emp}{\end{minipage}}
    \newcommand{\btb}{\begin{tabular}}
    \newcommand{\etb}{\end{tabular}}
    \newcommand{\barr}{\begin{array}}
    \newcommand{\earr}{\end{array}}
    \newcommand{\bit}{\begin{itemize}}
    \newcommand{\eit}{\end{itemize}}
    \newcommand{\ben}{\begin{enumerate}}
    \newcommand{\een}{\end{enumerate}}
    \newcommand{\bct}{\begin{center}}
    \newcommand{\ect}{\end{center}}
    \newcommand{\bfr}{\begin{flushright}}
    \newcommand{\efr}{\end{flushright}}
    \newcommand{\bea}{\begin{eqnarray*}}
    \newcommand{\eea}{\end{eqnarray*}}
    \newcommand{\bqo}{\begin{quote}}
    \newcommand{\eqo}{\end{quote}}
    \newcommand{\bdc}{\begin{description}}
    \newcommand{\edc}{\end{description}}
    \newcommand{\bdia}{\begin{CD}}
    \newcommand{\edia}{\end{CD}}

    \definecolor{light}{gray}{.3}

    \newcommand{\tr}{\mathrm{tr\,}}

    \newcommand{\ind}{\mathrm{ind\,}}
    
    \newcommand{\im}{\mathrm{im\,}}
    \newcommand{\res}{\mathrm{res\,}}

    \newcommand{\sr}[2]{{\,\stackrel{#1}{#2}\,}}
    \newcommand{\fra}[2]{{\,\frac{#1}{#2}\,}}

    \newcommand{\theorem}{\sn
           \bdc
           \item[{\sc Theorem.}] \em }
    \newcommand{\Theorem}[1]{\sn
           \bdc
           \item[{\sc Theorem {#1}.}]  \em }

    \newcommand{\Stop}{\edc \sn\rm} 

    \newcommand{\Lemma}[1]{\sn
           \bdc
           \item[{\sc Lemma {#1}.}] \em }
    
    \newcommand{\Proposition}[1]{\sn
           \bdc
           \item[{\sc Proposition {#1}.}] \em }

    \newcommand{\proof}{{\sc Proof.} \ }
    \newcommand{\remark}{\sn{\sc Remark.} \ }

\mmargin
\def\tr{{\mr{tr}}}
\def\dex{{\de^{(x)}}}
\def\ver{{\mr{ver}}}
\def\ab{{^{\mr{ab}}}}
\def\cn{{{\cal{N}}}}
\def\delg{{\De_g}}
\def\zpgs{{\zp[[G_S]]}}
\def\gsu{{G_S/U}}
\def\pm{pseudomeasure\ }
\def\pms{pseudomeasures\ }
\def\lk{{\la_K}}
\def\ll{{\la_L}}
\def\lg{{\la_g}}
\def\zphs{{\zp[[H_S]]}}
\def\ga{{g_\fa}}
\def\zek{{\ze_K(1-k,\dex)}}
\def\zeke{{\ze_K(1-k,\ve)}}
\def\zekeg{{\ze_K(1-k,\ve_g)}}
\def\muu{{m(U)}}
\def\pmu{{p^\muu}}
\def\cnl{{\cn_L}}
\def\mku{{m_K(U)}}
\def\pmku{{p^\mku}}
\def\mlv{{m_L(V)}}
\def\pmlv{{p^\mlv}}
\def\delgk{{\De_{g_K}}}
\def\delhl{{\De_{h_L}}}
\def\dekx{{\de_K^{(x)}}}
\def\dely{{\de_L^{(y)}}}
\def\cnk{{\cn_K}}
\def\lhl{{\la_{h_L}}}
\def\lgk{{\la_{g_K}}}
\def\abcd{{\bigl(\begin{smallmatrix}a&b\\ c&d \end{smallmatrix}\bigr)}}
\def\gan{{\Ga_{00}(\ff)}}
\def\sl{{\mr{SL}}}
\def\gafn{{\Ga_{00}(\ff\fo_L)}}

\begin{document}

\title{Congruences between abelian \pms}
 \author{Jürgen Ritter \ $\cdot$ \ Alfred Weiss \
 \thanks{We acknowledge financial support provided by NSERC and the University of Augsburg.}
 }
\date{\pht\today}

 \maketitle
 \bct\footnotesize{Dedicated to Professor Peter Roquette
 on his 80$^\mr{th}$ birthday}\ect

\bbn In this paper $K$ is a totally real number field (finite over
$\Q$), $p$ a fixed odd prime number, and $S$ a fixed finite set of
non-archimedian primes of $K$ containing all primes above $p$. Let
$K_S$ denote the maximal abelian extension of $K$ which is
unramified (at all non-archimedean primes) outside $S$ and set
$G_S=G(K_S/K)$. Serre's \pm $\lk=\la_{K,S}$ has the property that
$(1-g)\lk$ is in the completed group ring $\zpgs$ for all $g\in
G_S$ [Se2].
 \mn Let $L$ be a totally real Galois extension of $K$ of degree
 $p$ with group $\Si=G(L/K)$ and be such that all primes of $K$ which ramify in $L$
 are in $S$. Then $L_S\,,\,H_S=G(L_S/L)\,,\,\ll$
 are the corresponding objects over $L$ with respect to the set of
 primes of $L$ above $S$. Observe that $L_S$ is a Galois extension
 of $K$, with group $\fG$, hence $G_S$ is its maximal abelian
 factor group. And $H_S$ is a normal subgroup of $\fG$ of index
 $p$, so that $\Si=\fG/H_S$ acts on $H_S$ by conjugation. This
 situation induces the transfer map $\ver:G_S\to H_S$ by means of
 which we can compare $\lk$ and $\ll$.
 \theorem For $g_K\in G_S$ and
 $h_L=\ver(g_K)\in H_S$\,, $$\ver(\la_{g_K})\equiv\la_{h_L}\mod
 T\ ,$$
 where $\la_{g_K}=(1-g_K)\la_{K,S}\,,\,\la_{h_L}=(1-h_L)\la_{L,S}$\,, and
 where $T$ is the ideal in the ring $\zphs^\Si$ of $\Si$-fixed points
 of $\zphs$ consisting of all $\Si$-traces \ $\sum_{\si\in
 \Si}\al^\si\,,\,\al\in\zphs$\,.\Stop
 The proof follows from Deligne and Ribet [DR] by interpreting
 it on the Galois
 side as in [Se2]. Explicitly, the group $\ilim{}G_\ff$
 of [DR, p.230], with $\ff$ running through the integral ideals of $K$
 with all prime factors in $S$, is identified with our $G_S$, via class field
 theory.
 \bn The theorem implies the ``torsion congruences'' of [RW, \S3]
 in general. More precisely, let $L_\infty$ be the cyclotomic $\zp$-extension
 of $L$ and $\Ga_L=G(L_\infty/L)$. Choosing $g_K$ above to induce
 an automorphism of infinite order on $K_\infty$, the image
 of $1-h_L$ under $\zphs\to\zp[[\Ga_L]]$ is not in
 $p\zp[[\Ga_L]]$. Letting $\zphs_\bullet$ be the localization
 obtained by inverting the multiplicative set of elements of $\zphs$ whose
 image in $\zp[[\Ga_L]]$ is not in $p\zp[[\Ga_L]]$, the theorem
 reads $$\ver(\la_{K,S})\equiv\la_{L,S}\mod T_\bullet$$ with
 $T_\bullet$ the $\Si$-trace ideal in $\zphs^\Si_\bullet$\,. If
 $M$ is totally real and Galois over $K$ with $L_\infty\subset M\subset L_S$
 and
 $[M:L_\infty]$ finite, then the ``torsion congruence'' is
 obtained by specializing $\fG\to G(M/K)$. Moreover, for
 $p$-extensions $M/K$ the ``torsion congruences'' also imply the ``logarithmic congruences'' of [RW], so that
 we get a proof of many cases of the ``main conjecture'' of [RW2,FK]
 complementing the Heisenberg extensions of Kato \footnote{see
 [RW, `Added in proof']}.
 \mn Here is a short description of the individual sections to
 follow. In \S1 we write $\lgk$ as a limit element in
 $\ilim{U}\zp[G_S/U]/\pmu$ with $U$ open in $G_S$ and with
 certain integers $m(U)$. This allows us to study the claimed
 congruence on finite level, which is carried out in \S2. The next
 section is some preparation concerning Hilbert modular forms that
 we need for the proof of the theorem, which in \S4 is combined
 with the work of Deligne and Ribet to finish the proof. A final
 section briefly discusses a weaker version of the theorem when
 $p=2$.
 \Section{1}{Approximations to \pms}
 We review the construction of \pms [Se2] in a more explicit form
 that will be essential for our purposes. We first fix notation.
 \sn For a coset $x$ of an open subgroup $U$ of $G_S$ set
 $\dex(g)=1$ or 0 according to the cases $g\in x$ or $g\notin x$. Then, for integers
 $k\ge1$, define \ $\ze_K(1-k,\dex)=\ze_{K,S}(1-k,\dex)\in\Q$ \ to
 be the value at $1-k$ of the partial $\ze$-function for the set
 of integral ideals $\fa$ of $K$ prime to $S$ with Artin symbol
 $\ga$ in $x$ \footnote{so $\zek=\ze_S(x,1-k)$ in [Se2],  and
 $=L(1-k,\dex)$ in [DR], up to identification}. Note that the
 definition of $\zek$ extends  linearly to locally constant functions $\ve$
 on $G_S$ with values in a $\Q$-vector space and gives values
 $\ze_K(1-k,\ve)$ in that vector space.
 \mn Let $\cn=\cn_{K,p}:G_S\to\zp\mal$ be that continuous character whose
 value on $\ga$ for an integral ideal $\fa$ prime to $S$ is its
 absolute norm $\cn\fa$ \, \footnote{i.e., $\cn_p$ is the cyclotomic
 character, so $\cn_p(g)$ is determined by the action of $g$ on $p$-power roots of unity (see [Se2, (2.3)]}.
 For $g\in G_S\,,\,k\ge1$ and $\ve$ a locally constant
 $\qp$-valued function on $G_S$ we define, following [DR],
 $$\delg(1-k,\ve)=\zeke-\cn(g)^k\zekeg\in\qp\,,$$
 where $\ve_g(g')=\ve(gg')$ for $g'\in G_S$.
 \sn We can now state
 \Theorem{{\rm[(0.4) of [DR]]}} Let $\ve_1,\ve_2,\ldots$ be a
 finite sequence of locally constant functions $G_S\to\qp$ so that
 \ $\sum_{k\ge1}\ve_k(g')\cn(g')^{k-1}\in\zp$ for all $g'\in G_S$\,. Then
  $$\sum_{k\ge1}\delg(1-k,\ve_k)\in\zp\quad\mr{for\ all}\ g\in G_S\,.$$\Stop
 Call an open subgroup $U$ of $G_S$ {\em admissible}, if $\cn(U)\subset
 1+p\zp$, and define $\muu\ge1$ by
 $\cn(U)=1+p^\muu\zp$\,.
 \Lemma{1} If $U$ runs through the cofinal system of admissible open subgroups of $G_S$,
 then $\zpgs=\ilim{U}\zp[\gsu]/p^\muu\,\zp[\gsu]$\,.\Stop
 \proof The natural map
 $$\zpgs=\ilim{U}\zp[\gsu]\to\ilim{U}\zp[\gsu]/\pmu$$ is injective,
 since $m(U_1)\ge m(U_2)$ for $U_1\le U_2$ and since the $\muu$\,'s
 are unbounded. In order to show surjectivity, it is sufficient to
 find a linearly ordered cofinal family $\{U'\}$ of open subgroups,
 because then it follows that
 the image of $\ilim{U'}\zp[G_S/U']/p^{m(U')}$ is dense
 in the compact group
 $\ilim{U'}\zp[G_S/U']/p^{m(U')}$\,, by taking
 successive approximations which are compatible with the
 projections.
  Now, $G_S$ is finitely generated
 (over $\hat\Z$), as the inertia groups for the $\fp\in S$ are
 finitely generated and they together generate an open subgroup
 (the fixed field of which is the strict Hilbert class field of
 $K$). Thus $G_S$ is a homomorphic image of a finite product
 $\prod\hat\Z$ and hence the closed subgroup $(G_S)^{n!}$ has index
 dividing the finite order of $\prod(\hat\Z/n!)$ and so is
 open.
 \Proposition{2} For $g\in G_S$ there is a unique element
 $\lg\in\zpgs$, independent of $k$ \,
 \footnote{This allows us to take $k>2$ to avoid difficulties with $K=\Q$.}, whose image in
 $\zp[\gsu]/\pmu$ is $$\sum_{x\in\gsu}\delg(1-k,\dex)\cn
 (x)^{-k} x\mod\pmu\zp[G/U]$$ for all admissible $U$,
 where $\cn$ here also denotes the homomorphism $G_S/U\to(\zp/\pmu)\mal$ induced
 by our previous $\cn$. Moreover, if $\la$ is the \pm of
 [Se2], then $$(1-g)\la=\lg\,.$$ \Stop
 Note first that the displayed elements are well-defined by the
 definition of $\muu$ and that, varying $U$, they determine a
 limit element $\lg\in\zpgs$, since $\delg$ is a $\zp$-valued
 distribution (see [DR, (0.5)]).
 \sn We check that $\lg$ is independent of $k$. Fix $U$ and a
 coset $x$. Choose a (set) map $\eta:\gsu\to\zp\mal$ so that
 $\eta(g'U)\equiv\cn^{k-1}g'\mod\pmu$ for all $g'\in G_S$. Viewing
 $\eta$ as a locally constant function on $G_S$, then
 $$\delg(0,\dex\eta)\equiv\delg(1-k,\dex)\mod\pmu\,.$$
 To see this, apply Theorem [(0.4) of [DR]], repeated above, with
 $\ve_1=p^{-\muu}\dex\eta\,,\,\ve_k=-p^{-\muu}\dex$ ( and the other
 $\ve$\,'s zero). Hence, with $\ti x\in x$, $$\barr{l}\delg(1-k,\dex)\cn
 \ti x^{-k}\equiv\delg(0,\dex\eta)\eta(x)\me\cn \ti x\me\\
 =\delg(0,\eta(x)\me\dex\eta)\cn \ti x\me=\delg(0,\dex)\cn \ti x\me
 \mod\pmu\ .\earr$$
 We next check that our $\lg$ satisfies
 $$\group{\ve\cn^k,\lg}=\delg(1-k,\ve)$$ (compare [Se2, (3.6)]).
 As above, choose $\eta:\gsu\to\zp\mal$ so that now
 $\eta(yU)\equiv\cn^ky\mod\pmu$\,. Then, by [Se2, (1.1)],
 $$\barr{l} \group{\ve\cn^k,\lg}\equiv\group{\ve\eta,\lg}\equiv\sum_x\ve\eta(x)\delg(1-k,\dex)\cn
 x^{-k}\\
 \equiv\sum_x\ve(x)\delg(1-k,\dex)=\delg(1-k,\sum_x\ve(x)\dex)=\delg(1-k,\ve)\mod\pmu\
 .\earr$$
 By the argument following [Se2, (3.6)] it follows that $(1-g)\la$
 is equal to our $\lg$ for all $g\in G_S$.
 \Section{2}{Transfer}
 Let $L/K$ be as in the introduction. We decorate objects which
 depend on $L$ and are analogous to the ones of $K$ appropriately,
 e.g. $\cnl,m_L,\ldots$; in particular we have the notion of
 admissible open subgroups of $H_S$. Note that if $V$ is such an admissible open
 subgroup of $H_S$, then $\bigcap_{\si\in\Si}V^\si$ is also open
 and therefore the system of $\Si$-stable admissible open
 subgroups of $H_S$ is a cofinal system of open subgroups of
 $H_S$.
 \Lemma{3}\ben\item If $V$ is an admissible open
 subgroup of $H_S$ and $U$ is an admissible open
 subgroup of $G_S$ contained in $\ver\me(V)$, then
 $\mku\ge\mlv-1$\,.
 \item Let $y$ be a coset of a $\Si$-stable admissible open
 subgroup of $H_S$. If $h\in H_S$ is fixed by $\Si$, then
 $\De_h(1-k,\de_L^{(y^\si)})=\De_h(1-k,\dely)$\,, where
 $\De_h=\De_{L,h}$. In particular, $\lhl$ is fixed by $\Si$.\een\Stop
 The first assertion uses $\cnl(\ver(g))=\cnk (g)^p$ for $g\in G_S$.
 Now $U\le\ver\me(V)$ implies $\ver(U)\le V$, hence
 $\cnl(V)\supset\cnl(\ver(U))=\cnk(U)^p$, i.e.,
 $1+\pmlv\zp\supset(1+\pmku\zp)^p=1+p^{\mku+1}\zp$\,. Thus
 $\mku+1\ge\mlv$.
 \sn For the second assertion it suffices to show that
 $\ze_L(1-k,\de_L^{(y^\si)})=\ze_L(1-k,\de_L^{(y)})$ for all $y$,
 because $(\de_L^{(y^\si)})_h=\de_L^{(h\me y^\si)}=\de_L^{((h\me
 y)^\si)}$ and $\de_L^{(h\me y)}=(\de_L^{(y)})_h$\,. Now
 view $\dely$ as a complex valued
 function on $H_S/V$ and write it as a $\C$-linear combination of
 the (abelian) characters $\chi$ of $H_S/V$. It suffices
 to check whether $\ze_L(1-k,\chi)=\ze_L(1-k,\chi^\si)$\,, with
 $\chi^\si(h)=\chi(h^{\si\me})=\chi(\si h \si\me)$. But
 this follows from the
 compatibility of the Artin $L$-functions with induction, because
 $\ind_{H_S/V}^{\fG/V}\chi=\ind_{H_S/V}^{\fG/V}\chi^\si$\,.
 \sn This finishes the proof.
 \mn Let $N$ be the kernel of $\ver:G_S\to H_S$. A $\Si$-stable admissible open subgroup
 $V$ of $H_S$ gives rise to the transfer map $\gsu\to H_S/V$
 whenever $U\le\ver\me(V)$. These transfer maps combined yield the
 right hand map in the commutative square
 $$\barr{cccc}\zpgs&\to&\ilim{U\ge N}\zp[\gsu]/\pmku&\\
 \daz{\ver}&&\da&\\
 \zphs&\sr{\simeq}{\lto}&\ilim{V,\Si\!-\!\mr{stable}}\zp[H_S/V]/p^{\mlv-1}&,\earr$$
 explicitly sending $(x_U)_U$ to $(y_V)_V$ where $y_V$ is the image
 of $x_U$ under
 $\zp[\gsu]/\pmku\sr{\ver}{\lto}$
 $\zl[H_S/V]/\pmku\to\zp[H_S/V]/p^{\mlv-1}$ whenever
 $U\le\ver\me(V)$.
 The bottom arrow is an isomorphism by the proof of Lemma 1.
 \mn We recall that a locally constant function $\ve_L$ on $H_S$
 is {\em even}, if $\ve_L(c_wh)=\ve_L(h)$ for all $h\in H_S$ and all
 ``Frobenius elements'' $c_w$ at the archimedean primes $w$ of $L$ (so
 $c_w\in H_S$ is the restriction of complex conjugation with
 respect to an embedding $L_S\hookrightarrow\C$ inducing $w$ on
 $L$).
 \sn Set $\Z_{(p)}=\Q\cap\zp$.
 \Proposition{4} A sufficient condition
 for the Theorem in the introduction to hold is the
 following\,:
 $$\delhl(1-k,\ve_L)\equiv\delgk(1-pk,\ve_L\circ\ver)\mod p\zp$$
 for all even locally constant $\Z_{(p)}$-valued functions $\ve_L$ on $H_S$
 satisfying $\ve_L^\si=\ve_L\ (\forall\,\si\in\Si)$ with
 $\ve_L^\si(h)=\ve_L(h^{\si\me})$.\Stop
 \proof Look at the coordinates of $$\lhl\
 \ \mr{and}\ \ \ver(\lgk)\ \ \mr{in}\ \
 \zp[H_S/V]/p^{\mlv-1}$$ for a
 $\Si$-stable admissible open subgroup $V\le H_S$ containing the group $C$ generated by all
 elements $c_w$. Note that $\ver(\lgk)$ is then the image
 under `ver' of the $U$-coordinate of $\lgk$, where $U=\ver\me(V)\le G_S$ contains $N$.
  These coordinates are the images of
 \bit\item[(i)]
 $\sum_{y\in H_S/V}\delhl(1-k,\dely)\cnl(y)^{-k}y$\,,\nz
 respectively
 \item[(ii)]
 $\sum_{x\in\gsu}\delgk(1-pk,\dekx)\cnk(x)^{-pk}\ver(x)$
 \eit  in $(\zp[H_S/V]/p^{\mlv-1})^\Si$ by Proposition 2 (recall that it asserts
 independence of $\lg$ from $k$).
 \sn We show that the sums in (i),(ii) are congruent modulo
 $T(V)$, where $T(V)$ is the $\Si$-trace ideal in
 $(\zp[H_S/V]/p^{\mlv-1})^\Si$, by
 distinguishing two cases\,:
 \ben \item {\em $y$ is fixed by $\Si$.} \ Then $\dely$ is an $\ve_L$ as
 appearing in the proposition and so\linebreak
 $\delhl(1-k,\dely)\equiv\delgk(1-pk,\dely\circ\ver)\mod p\,.$
 Now, if $y=\ver(x)$, then, because  $\ver:U/N\to V$ is an isomorphism, $x$ is uniquely determined by
 $y$ and
 $\cnl(y)^{-k}=\cnl(\ver(x))^{-k}=\cnk(x)^{-pk}$\,.
 Moreover, $\dely\circ\ver=\dekx$.
 Hence the corresponding summands in (i) and (ii) cancel out
 modulo $T(V)$, since $p\al$ is a $\Si$-trace whenever $\al$ is
 $\Si$-invariant.
 However, if $y\notin\im(\ver)$, then
  $\dely\circ\ver=0$, hence the
 $y$-summand vanishes modulo $T(V)$.
 \item {\em $y$ is not fixed by $\Si$.} \ By 2.~of Lemma 3,
 $\delhl(1-k,\dely)=\delhl(1-k,\de_L^{(y^\si)})$, whence the
 $\Si$-orbit of $y$ yields the sum
 $\delhl(1-k,\dely)\cnl(y)^{-k}\sum_{\si\in\Si}y^\si$\, which is
 in $T(V)$.\een
 Now subtracting type (ii) sums from type (i) sums for all $\Si$-stable admissible open
 $V\ge C$ gives a compatible system of elements in $\ilim{V\ge
 C}T(V)\subset\ilim{V\ge C}\zp[H_S/V]/p^{\mlv-1}$. Set $H_S^+=H_S/C$; so $H_S^+=G(L_S^+/L)$ where $L_S^+$ is the maximal
 totally real subfield of $L_S$.  Since
 $T(V_1)\to T(V)$ is surjective whenever $V_1\le V$, we get a
 limit $s^+\in T^+\subset \zp[[H_S^+]]$.
  Thus the  proposition follows from
 \Lemma{5} Suppose that
 $s\df\lhl-\ver(\lgk)\in\zphs^\Si$ has image $s^+$ under $\zphs\to\zp[[H_S^+]]$ in the
 $\Si$-trace ideal $T^+$ in $\zp[[H_S^+]]^\Si$. Then $s\in T$.
 \Stop
 \proof We know, from [Se2, (3.12)], that the Frobenius elements
$c_v\in G_S$ for the
 real primes $v$ of $K$ satisfy $c_v^2=1\,,\,c_v\lgk=\lgk$\,, and
 that they generate the kernel of $G_S\to G_S^+$. Put
 $c_K=\prod_v(1+c_v)$.
 \sn The analogous properties hold for the $c_w$ for the real
 primes $w$ of $L$, and we can form $c_L$. Moreover,
 $$c_L\equiv\ver(c_K)\mod T\,.$$
 To see this, expand $c_L$ in a sum of products of $c_w$\,'s and
 consider the $\Si$-action on the summands. The sum of each orbit
 of length $p$ is in $T$ and the products fixed by $\Si$ add up
 to $\ver(c_K)$, because $\ver(c_v)=\prod_{w|v}c_w$ for every $v$.
 \sn Now $s^+\in T^+$ and the surjectivity of $T\to T^+$ mean that
 $s^+=t^+$ for some $t\in T$, hence $s-t$ is in the kernel of
 $\zphs\to\zp[[H_S^+]]$ which is generated by all $1-c_w$ as a $\zphs$-module. Then
 $c_L(s-t)=0$, implying
 $c_Ls\in T$ because $c_L\in\zphs^\Si$ and $T$ is an ideal of
 $\zphs^\Si$.
 \sn Moreover $c_L\lhl=2^{[L:\Q]}\lhl$ and
 $$c_L\ver(\lgk)\equiv\ver(c_K)\ver(\lgk)=\ver(c_K\lgk)=2^{[K:\Q]}\ver(\lgk)\mod
 T\ .$$ Since $2^p\equiv2\mod p$, it follows that
 $2^{[K:\Q]}s\equiv c_L(\lhl-\ver(\lgk))=c_Ls\equiv0\mod T$\,,
 from which the lemma follows as $p\in T$ is odd.
 \Section{3}{$q$-expansions}
 Let $[K:\Q]=r$, let $\ff$ be an integral
 ideal with all prime factors in $S$, and denote the usual Hilbert upper half plane
 associated to $K$ by $\fH=\{\tau\in
 K\ot\C:\Im(\tau)\gg0\}$.
 \sn If $k$ is an even positive integer, we define, as usual, the action of
 $\mr{GL}^+(2,K\ot\R)$ of matrices with totally positive
 determinant on functions $F:\fH\to\C$ by
 $$(F_{|k}\abcd)(\tau)=\cn(ad-bc)^{k/2}\cn(c\tau+d)^{-k}F(\fra{a\tau+b}{c\tau+d})\,,$$
 with $\cn:K\ot\C\to\C$ denoting the norm.
 \sn Set
 $$\Ga_{00}(\ff)=\{\abcd\in\mr{SL}(2,K):a,d\in1+\ff\,,\,b\in\fD\me\,,\,c\in\ff\fD\}$$
 where $\fD$ is the different of $K$. A Hilbert modular form $F$ of
 weight $k$ on $\gan$ is a holomorphic function $\fH\to\C$ \,
 \footnote{and holomorphic at infinity, if $K=\Q$}
 satisfying $F_{|k}M=F$ for all $M\in\gan$. Denote the space of these by
 $M_k(\gan,\C)$ (see [DR, (5.7)]). Such an $F$ can be
 expanded as a Fourier series $$c(0)+\sum_{\substack{\mu\in\fo_K\\
 \mu\gg0}}c(\mu)q^\mu\quad\mr{with}\quad q^\mu=e^{2\pi
 i\tr(\mu\tau)}\quad \footnote{$\fo_K$ is the ring of integers in $K$; from now on $\mu$ will always be in $\fo_K$}\,,$$
 called the standard $q$-expansion of $F$, i.e., the
 $q$-expansion at the cusp $\infty=\fra10$.
 \Lemma{6} Let $\be\in\fo_K$ be totally positive with $\ff\subset\be\fo_K$. There is a Hecke
 operator $U_\be$ on $M_k(\gan,\C)$ so that, if $F\in
 M_k(\gan,\C)$ has standard $q$-expansion as above, then
 $F_{|k}U_\be\in M_k(\gan,\C)$ has standard $q$-expansion
 $c(0)+\sum_{\mu\gg0}c(\be\mu)q^\mu$\,.\Stop
 Following [AL, \S\S2,3] for the proof, let $B=\bigl(\begin{smallmatrix}\be&0\\
 0&1\end{smallmatrix}\bigr)$ and set $\Om=B\gan B\me\cap\gan$\,.
 The matrices $S_\xi=\bigl(\begin{smallmatrix}1&\xi\\
 0&1\end{smallmatrix}\bigr)$, with $\xi$ running through a set of
 coset representatives of $\be\fD\me$ in $\fD\me$, satisfy
 $\gan=\dot\bigcup_\xi\Om S_\xi$\,, because $\ff\subset\be\fo_K$.
 \sn Define $U_\be$ on $M_k(\gan,\C$) by
 $$F_{|k}U_\be=\cn(\be)^{\fra{k}{2}-1}\sum_\xi F_{|k}B\me
 S_\xi\,,\ \mr{with}\ \cn\ \mr{as\ above\ in\ this\ section.}$$
 Then $F_{|k}B\me$ is modular on $B\gan B\me$, hence on $\Om$. The
 usual averaging argument then shows that $F_{|k}U_\be$ is modular
 on $\gan$. Now,
 $$\barr{l} (F_{|k}U_\be)(\tau)=\cn(\be)\me\sum_\xi
 F(\be\me\tau+\be\me\xi)=\cn(\be)\me\sum_\xi\Big(c(0)+\sum_{\mu\gg0}c(\mu)e^{2\pi
 i\tr(\mu\fra{\tau+\xi}{\be})}\Big)\\
 =\cn(\be)\me[\fD\me:\be\fD\me]c(0)+\sum_{\mu\gg0}c(\mu)\Big(\cn(\be)\me\sum_\xi e^{2\pi
 i\tr_{K/\Q}(\mu\xi/\be)}\Big)e^{2\pi
 i\tr(\mu\tau/\be)}\\=c(0)+\sum_{\substack{\mu\gg0\\
 \mu\in\be\fo_K}}c(\mu)e^{2\pi i\tr(\mu\tau/\be)}\earr$$
 as $\cn(\be)=[\fD\me:\be\fD\me]$ and as $\xi\mapsto e^{2\pi
 i\tr_{K/\Q}(\mu\xi/\be)}$ is a character on $\fD\me/\be\fD\me$.
 Since $\be$ is totally positive, the proof of the lemma is
 complete.
 \bn We next discuss restriction of Hilbert modular forms from $L$
 to $K$.
 \sn The containment $K\subset L$ induces natural maps
 $\fH_K\sr{\ast}{\to}\fH_L$ and
 $\sl(2,K\ot\R)\sr{\ast}{\to}\sl(2,L\ot\R)$\,. For a
 holomorphic $F:\fH_L\to\C$ define the restriction $\res F:\fH_k\to\C$ of $F$ to be
 the holomorphic function satisfying \ $(\res
 F)(\tau)=F(\tau^\ast)$\,. Then
 \nf{$\star$}{(\res F)_{|pk}M=\res(F_{|k}M^\ast)\ \mr{for}\
 M\in\sl(2,K\ot\R)\,.}
 \noi The $q$-expansion at a cusp determined by a finite id\`ele
 $\al\in\hat K\mal$ is discussed in [DR, bottom of p.229 and (5.8)].
 \Lemma{7}\quad Let $F\in M_k(\gafn,\C)$ and let
 \ $c(0)+\sum_{\substack{\nu\gg0\\
 \nu\in\fo_L}}c(\nu)q_L^\nu$ \ be its standard $q$-expansion (with $q_L^\nu=e^{2\pi i\tr_L(\nu\tau)}$).
 Let $\al\in\hat K\mal$\,.
 Then \ben\item\quad  $\res F\in M_{pk}(\gan,\C)$ \ has
 standard \ $q$-expansion\quad $c(0)+\sum_{\substack{\mu\gg0\\ \mu\in\fo_K}}c_\ast(\mu)q_K^\mu$\hsp{4} \linebreak with
 $c_\ast(\mu)=\sum_{\substack{\nu\gg0\,,\,\nu\in\fo_L\\
 \tr_{L/K}(\nu)=\mu}}c(\nu)$ (and $q_K^\mu=e^{2\pi i\tr_K(\mu\tau)}$)\,,
 \item\quad the constant term of \ $\res F$ \ at the cusp determined by
 $\al$ equals the constant\hsp{5}\linebreak term of $F$ at the cusp determined by
 $\al^\ast\in\hat L\mal$\,.\een\Stop
 Assertion 1.~follows from observing that
 $\tr_L(\nu\tau^\ast)=\tr_K(\tr_{L/K}(\nu)\tau)$ for $\nu\in
 L\,,\,\tau\in K\ot\C$, and substituting this into the definition.
 \sn For 2., the constant terms in
 question are those of
 $$(\res F)_\al=(\res F)_{|pk}\bigl(\begin{smallmatrix}\al&0\\
 0&\al\me\end{smallmatrix}\bigr)\quad\mr{and}\quad F_{\al^\ast}
 =F_{|k}\bigl(\begin{smallmatrix}\al&0\\
 0&\al\me\end{smallmatrix}\bigr)^\ast\,,$$
 respectively, by [DR, p.229].
 By 1., $F_{\al^\ast}$ and $\res F_{\al^\ast}$ have the same
 constant term in their respective standard $q$-expansion, so it
 suffices to show $(\res F)_\al=\res F_{\al^\ast}$. For that,
 decompose $M=\bigl(\begin{smallmatrix}\al&0\\
 0&\al\me\end{smallmatrix}\bigr)\in\sl(2,\hat K)$ as
 $M=M_1M_2$ according to $\sl(2,\hat
 K)=\widehat{\gan}\cdot\sl(2,K)$\,, hence
 $M^\ast=M_1^\ast M_2^\ast$ according to $\sl(2,\hat
 L)=\widehat{\gafn}\cdot\sl(2,L)$. Then
 $$(\res F)_\al=(\res F)_{|pk} M =(\res
 F)_{|pk}M_2\sr{(\star)}{=}\res(F_{|k}M_2^\ast)=\res(F_{|k}M^\ast)=\res
 F_{\al^\ast}\,,$$ with equation $\sr{(\star)}{=}$ referring to
 the formula displayed prior to Lemma 7.
 \Section{4}{Proof of the main result}
 We use the notation of the previous section, except that we now also use $\cn$ for the norm map $K\to\Q$ and any norm map
 derived from it \footnote{hence consistent with our usage in \S3}, as in [DR,\S2].
 \sn
 We
 attach an Eisenstein series of every even weight $k$ to even locally constant $\C$-valued
 functions $\ve$ via [DR, (6.1)].
 \Proposition{8}Let $\ve$ be an even locally constant $\C$-valued function on $G_S$.
 \ben\item There is an integral ideal $\ff$ in $K$ with all its
 prime factors in $S$ and a modular form
 $G_{k,\ve}\in M_k(\gan,\C)$ with standard $q$-expansion
 $$2^{-r}\ze_K(1-k,\ve)+\sum_{\substack{\mu\gg0\\ \mu\in\fo_K}}\Big(\sum_{\substack{\mu\in\fa\subset\fo_K\\
 \fa\,\mr{prime\,to}\, S}}\ve(\fa)\cn(\fa)^{k-1}\Big)q^\mu$$
 where $\ve(\fa)=\ve(g_\fa)$ with $g_\fa\in G_S$ the
 Artin symbol of $\fa$.
 \item Its $q$-expansion at the cusp determined by
 $\al\in\hat K\mal$ has constant term
 $$\cn((\al))^k2^{-r}\ze_K(1-k,\ve_{a})\,,$$ where $(\al)$ is
 the ideal generated by $\al$ and $a\in G_S$ is the image of $\al$ under the map
 \nf{2a}{\hat K\mal\sr{j}{\lto}G= G(K\ab/K)\onto G_S}
 with $j$ taken from [DR, (2.22)] and the identification
 $G=G(K\ab/K)$ as in [DR, p.240], via the Artin symbol on integral
 ideals prime to $\ff$.
 \item $\cn((\al))=\cn(\al_p)\cdot\cn_p(a)$ where $\al_p\in K\ot\qp$ is the $p$-component
 of $\al\in\hat K\mal$ and $\cn_p(a)=\cn_{K,p}(a)$\,, as in \S1. \een\Stop
 For 1. choose an open subgroup $U$ of $G_S$ so that $\ve$ is
 constant on each coset of $G_S/U$. Let $\ff$ be an integral ideal
 which is a multiple of the conductor of the field fixed by $U$
 acting on $K_S$ and with all its prime factors in $S$. Then the Artin symbol maps
 the strict ideal class group $G_\ff$ onto $G_S/U$. Viewing
 $G_\ff$ as the group of invertible elements of $A_\ff$, as in
 [DR, (2.6)], makes $\ve$ a map on $G=\ilim{\sss{\ff'\subset
 \ff}}G_{\ff'}$\,. Finally extend $\ve$ to $I$ by zero
 \footnote{for the definition of $I$ see [DR, \S2]}. In
 particular, if $\fa$ is an integral ideal prime to $S$, then
 $\ve(\fa)=\ve(g_\fa)$. Moreover, by [DR, (2.3) and (2.4)],
 $\ve(\fa)=0$ for every (fractional) ideal $\fa$ of $K$ which is
 not integral and prime to $S$.
 \sn Now, with this $\ve$, [DR, (6.2)] gives the standard
 $q$-expansion of $G_{k,\ve}$\,:
 $$2^{-r}\ze_K(1-k,\ve)+\sum_{\substack{\mu\gg0\\ \mu\in\fo_K}}(\sum_{\fx\subset\fo_K}\ve(\mu\fx\me)\cn(\mu\fx\me)^{k-1})q^\mu\,,$$
 where we have chosen the ideal $\fB$ of [DR] to be $\fo_K$. Set $\fa=\mu\fx\me$, so
 $\mu\in\fa$, and we may assume that $\fa$ is integral and prime
 to $S$, because otherwise $\ve$ will be zero on $\fa$. Thus the
 above $\mu$\,th coefficient is turned into \
 $\sum_{\substack{\mu\in\fa\subset\fo_K\\
 \fa\,\mr{prime\,to}\,S}}\ve(\fa)\cn(\fa)^{k-1}\,.$
 \mn For 2., [DR, (6.2)] shows that
 $\cn((\al))^k2^{-r}\ze_K(1-k,\ve_c)$ is the constant term of the
 $q$-expansion at the cusp determined by $\al\in \hat K\mal$, with
 $c=j(\al)$. Our extension of $\ve$ to $G$ has been such that, for
 $g\in G$, $\ve_c(g)=\ve(cg)=\ve(\ol c\ol g)=\ve_{\ol c}(\ol g)$
 with $\ol c,\ol g$ the images of $c,g$ in $\ilim{{\sss{\ff'\subset\ff}},
 \mr{in}\,S}G_{\ff'}$\,, where `\,$\ff'\,\mr{in}\,S$\,' means that every
 prime factor of $\ff'$ is in $S$. Hence the commutative square
 $$\barr{ccc}G&\lto&\ilim{{\sss{\ff'\subset\ff}},
 \mr{in}\,S}G_{\ff'}\\ \pl&&\pl\\ G(K\ab/K)&\lto&G_S\earr$$ shows
 $\ol c=a$, up to identification.
 \mn For 3., we get from [DR, (2.12),(2.16)] that the norm of
 $c=j(\al)$ is $\cn(\al)\me\cn((\al))$. Thus the $p$-component of
 $\cn(c)\in\hat \Z\mal$ in $\zp\mal$ is $\cn(\al_p)\me\cn((\al))$
 since $\cn((\al))\in\Q\mal$. On the other hand, the $p$-component
 of $\cn(c)$ is $\cn_p(a)$ by the commutative diagram
 $$\barr{ccccc}c\in
 G&\lto&\ilim{n}G_{p^n\fo_K}&\sr{\mr{N}_{K/\Q}}{\lto}&\ilim{n}G_{p^n\Z}\\
 \mapsdto\pht{\in}\sda&&&&\daz{\simeq}\\
 a\in
 G_S&\lto&\ilim{n}G(\Q(\mu_{p^n})/\Q)&\sr{\simeq}{\lto}&\ilim{n}(\Z/p^n)\mal\earr$$
 with the left map as in $(\sharp)$ and $\mu_{p^n}$ the $p^n$\,th roots of unity. Here, the map
 $G\to\ilim{n}(\Z/p^n)\mal=\zp\mal$ around the top row takes $c$
 to the $p$-component of $\cn(c)$, which thus is $\cn_p(a)$.
 \sn The proof of the proposition is complete.
 \Lemma{9}\footnote{compare [Ty]} Let $k$ be an even positive integer and $\ve_L$ an
 even locally constant $\Z_{(p)}$-valued function on $H_S$. There
 is an integral ideal $\ff\subset p\fo_K$ with all prime factors
 in $S$, so that
 $$E=(\res {G_{k,\ve_L}})_{|pk}U_p-G_{pk,\ve_L\circ\ver}\quad\mr{is\
 in}\quad M_{pk}(\gan,\C)\quad\footnote{recall that
 $U_p=\bigl(\begin{smallmatrix}1&p\\
 0&1\end{smallmatrix}\bigr)$}\,.$$
 If $\ve_L^\si=\ve_L$ for all $\si\in\Si$, then the constant term of the standard $q$-expansion of $E$ is
 $$2^{-pr}\ze_L(1-k,\ve_L)-2^{-r}\ze_K(1-pk,\ve_L\circ\ver)$$ and
 all non-constant coefficients are in $p\Z_{(p)}$.\Stop
 Choose an $\ff\subset p\fo_K$ by Proposition 8 so that $G_{pk,\ve_L\circ\ver}\in
 M_{pk}(\gan,\C)$ and $G_{k,\ve_L}\in M_k(\gafn,\C)$\,.
 Using Lemmas 6 and 7, the standard $q$-expansion of
 $G_{pk,\ve_L\circ\ver}$ is $$2^{-r}\ze_K(1-pk,\ve_L\circ\ver)+\sum_{\substack{\mu\gg0\\ \mu\in\fo_K}}\Big(\sum_{\substack{\mu\in\fa\subset\fo_K\\
 \fa\,\mr{prime\,to}\,S}}\ve_L(\fa\fo_L)\cn_K(\fa)^{pk-1}\Big)q_K^\mu$$
 because
 $(\ve_L\circ\ver)(\fa)=(\ve_L\circ\ver)(g_\fa)=\ve_L(\ver(g_\fa))=\ve_L(\fa\fo_L)$ (see [Se1, VII,8])\,,
 and that of
 $(\res{G_{k,\ve_L}})_{|pk}U_p$ is
 $$2^{-pr}\ze_L(1-k,\ve_L)+\sum_{\substack{\mu\gg0\\ \mu\in\fo_K}}\Big(\sum_{\substack{(\fb,\nu)\,\mr{so}\,
 \nu\in\fb\subset\fo_L\,,\,\nu\gg0\\
 \fb\,\mr{prime\,to}\,S\,,\,\tr_{L/K}(\nu)=p\mu}}\ve_L(\fb)\cn_L(\fb)^{k-1}\Big)q_K^\mu\
 .$$
 Hence, the $\mu$\,th coefficient of $E$ is
 $$\sum_{(\fb,\nu)}\ve_L(\fb)\cn_L(\fb)^{k-1}-\sum_\fa\ve_L(\fa\fo_L)\cn_K(\fa)^{pk-1}$$
 with $(\fb,\nu)$ so that $\nu\gg0,\nu\in\fb\subset\fo_L,\fb$ prime to $S$, $\tr_{L/K}(\nu)=p\mu$ and $\fa\subset\fo_K$
 prime to $S$.
 \sn The group $\Si$ acts on the pairs $(\fb,\nu)$ by
 $(\fb,\nu)^\si=(\fb^\si,\nu^\si)$. If $\Si$ moves $(\fb,\nu)$,
 then the orbit sum
 $\sum_\si\ve_L(\fb^\si)\cn_L(\fb^\si)^{k-1}=p\ve_L(\fb)\cn_L(\fb)^{k-1}$
 because $\ve_L(\fb^\si)=\ve_L^{\si\me}(\fb)=\ve_L(\fb)$ and
 $\cn_L(\fb^\si)=\cn_L(\fb)$.
 \sn However, if $\Si$ fixes $(\fb,\nu)$, then $\nu\in
 K\,,\,\tr_{L/K}(\nu)=p\mu,$ so $\nu=\mu$, and $\fb^\si=\fb$, so
 $\fb=\fa\fo_L$ for a unique integral ideal $\fa$ of $K$ prime to $S$,
 since $S$ contains all primes which are ramified in $L$. Thus
 $(\fb,\nu)=(\fa\fo_L,\mu)$.
 \sn The above claim on $E$ now follows from
 $$\ve_L(\fb)\cn_L(\fb)^{k-1}=\ve_L(\fa\fo_L)\cn_L(\fa\fo_L)^{k-1}=\ve_L(\fa\fo_L)\cn_K(\fa)^{p(k-1)}\\
 \equiv\ve_L(\fa\fo_L)\cn_K(\fa)^{pk-1}\mod p\ ,
 $$ by $\cn_K(\fa)^{p-1}\equiv1\mod p$\,.
 \bn
 We finally turn to the {\sc Proof} of the {\sc Theorem} stated in the
 introduction. We check the sufficient conditions for every
 $\ve_L$ as in Proposition 4. These are the $\ve_L$ appearing in
 Lemma 9. With $E$ as in Lemma 9 and $\al\in \hat K\mal$, let
 $E_\al$ be the $q$-expansion of $E$ at the cusp determined by
 $\al$ and let $E(\al)=\cnk(\al_p)^{-pk}E_\al$.
 \mn Since, by [DR, (2.23)], the map $j$ in $(2a)$ is surjective, there is an
 id\`ele $\gamma\in \hat K\mal$ which maps to $g_K\in G_S$ by
 $(2a)$.
 According to Lemma 9, $E(1)=E_1$ has non-constant coefficients in $p\Z_{(p)}$.
 Then, by [DR, (0.3) and {\em Variant: Forms on $\gan$} at the end of \S5],
 $E(1)-E(\gamma)$ has constant coefficient in $p\zp$. This
 coefficient is, by Lemmas 6,7,9 and Proposition 8,
 $$\barr{l}
 2^{-pr}\ze_L(1-k,\ve_L)-2^{-r}\ze_K(1-pk,\ve_L\circ\ver)\ -\\ \hsp{4}\cn_K(\gamma_p)^{-pk}\cn_K((\gamma))^{pk}
 \Big[2^{-pr}\ze_L(1-k,(\ve_L)_{h_L})-2^{-r}\ze_K(1-pk,(\ve_L\circ\ver)_{g_K})\Big]\\
 =2^{-pr}\Big[\ze_L(1-k,\ve_L)-\cn_K(g_K)^{pk}\ze_L(1-k,(\ve_L)_{h_L})\Big]\ -\\ \hsp{4}2^{-r}\Big[\ze_K(1-pk,\ve_L\circ\ver)
 -\cn_K(g_K)^{pk}\ze_K(1-pk,(\ve_L\circ\ver)_{g_K})\Big]\\
 =2^{-pr}\De_{h_L}(1-k,\ve_L)-2^{-r}\De_{g_K}(1-pk,\ve_L\circ\ver)\ \equiv\\
 \hsp{4}
 2^{-r}\Big(\De_{h_L}(1-k,\ve_L)-\De_{g_K}(1-pk,\ve_L\circ\ver)\Big)\mod
 p\earr$$
 where we have used that $\gamma^\ast\in\hat L\mal$ has image $\ver(g_K)=h_L$ under the map $(2a)_L$ as well as
 $\cn_K(g_K)^p=\cn_L(h_L)\ .$
 \sn Thus, Proposition 4 finishes the proof.
 \Section{5}{About $p=2$}
For $p=2$ the theorem needs to be reformulated because of the
``extra'' 2-adic divisibilities of [DR]. In view of Lemma 8, we
define $$\ti\ze_{K,S}(1-k,\ve)=2^{-r}\ze_{K,S}(1-k,\ve)\,, $$
whence \ $\ti\De_{g_K}(1-k,\ve)=2^{-r}\De_{g_K}(1-k,\ve)$ \ takes
values in $\Z_2$ for $\Z_2$-valued $\ve$, since an admissible
subgroup never admits conductor (1) (see [Ri, \S3]). Hence
$\ti\la_{g_K}=2^{-r}\lgk$ is in $\Z_2[[G_S]]$ (by e.g.~Proposition
2). Following the proof of the theorem now shows that $$\mr{the\
image\ of} \ \, \ver(\ti\la_{g_K})-\ti\la_{h_L} \ \,\mr{under} \
\, \Z_2[[H_S]]\to\Z_2[[H_S^+]] \ \,\mr{is\ in}\ \ \, T^+\ ,$$ in
the notation of Lemma 5. But the proof of Lemma 5 does not work
anymore. One imagines that the methods of [DR], which gave the
extra 2-adic divisibilities in the first place, would also sharpen
the conclusion displayed above.
 \remark Actually, we can do the same modification for odd $p$.
 The equivariant ``main conjecture'' of [RW2] is unaffected
 because $[{\cal{Q}}G_\infty,2]$ is then in the kernel of
 $\partial:K_1({\cal{Q}}(\zp[[G_\infty]]))\to K_0T(\zp[[G_\infty]])$
 (see equation $(\Im)$ on p.\,550 of [RW2]).
  \bbn{\large {\sc References}}
 \small
 \bn
 \btb{rp{13cm}}
 \,[AL]   & Atkin, A.O.L., Lehner, J., {\em Hecke Operators on
          $\Ga_0(m)$.} Math.~Ann.~{\bf 185} (1970), 134-160\\
 \,[DR]   & Deligne, P.~and Ribet, K., {\em Values of abelian
            $L$-functions at negative integers
            over totally real fields.} Invent.~Math. {\bf 59} (1980), 227-286\\
 \,[FK]   & Fukaya, T., Kato, K., {\em A formulation of conjectures on $p$-adic zeta
            functions in noncommutative Iwasawa theory.} Proceedings of
            the St.~Petersburg Mathematical Society, vol.~XII (ed.
            N.N.~Uraltseva), AMS Translations -- Series 2, {\bf 219} (2006),
            1-86\\
 \,[Ri]    & Ribet, K., {\em Report on $p$-adic $L$-functions over
            totally real fields.} Ast\`erisque {\bf 61} (1979), 177-192\\
 \,[RW]    & Ritter, J.~and Weiss, A., {\em Non-abelian pseudomeasures
           and congruences between abelian Iwasawa $L$-functions.}
           To appear in Pure and Applied Mathematics Quarterly\\
 \,[RW2]   & ---------------------, {\em Toward equivariant Iwasawa theory, II.}
           Indagationes Mathematicae {\bf 15} (2004), 549-572 \\
 \,[Se1]   & Serre, J.-P., {\em Corps locaux.} Hermann, Paris
            (1968)\\
 \,[Se2]   & ---------, {\em Sur le r\'esidu de la fonction
           z\^{e}ta $p$-adique d'un corps de nombres.}
           C.R.Acad.Sci.~Paris {\bf 287} (1978), s\'erie A, 183-188\\
 \,[Ty]  & Taylor, M.J., {\em Galois module type congruences for
            values of $L$-functions.} J.~LMS {\bf
            24} (1981), 441-448\\
 \etb

 \bn {\footnotesize \bct Institut für Mathematik $\cdot$
 Universität Augsburg $\cdot$ 86135 Augsburg $\cdot$ Germany \\
 Department of Mathematics $\cdot$ University of Alberta $\cdot$
 Edmonton, AB $\cdot$ Canada T6G 2G1   \ect

\end{document}